\documentclass{amsart}

\usepackage{amsmath,amssymb,amsfonts,enumerate,amsthm, amscd}

\newcommand{\pd}{\mbox{pd}\,}
\newcommand{\id}{\mbox{id}\,}
\newcommand{\fd}{\mbox{fd}\,}

\newcommand{\wdim}{\mbox{wdim}\,}
\newcommand{\gldim}{\mbox{gldim}\,}

\newtheorem{theorem}{Theorem}[section]
\newtheorem{lemma}[theorem]{Lemma}
\newtheorem{proposition}[theorem]{Proposition}
\newtheorem{corollary}[theorem]{Corollary}

\theoremstyle{definition}

\theoremstyle{remark}
\newtheorem{remark}[theorem]{Remark}

\newtheorem{discussion}[theorem]{Discussion}

\theoremstyle{Definition and Notation}

\begin{document}
\bibliographystyle{amsplain}

\title[On Gorenstein global dimension in trivial ring extensions]{On Gorenstein global dimension in trivial ring extensions}

\author{Najib Mahdou}
\address{Najib Mahdou\\Department of Mathematics, Faculty of Science and Technology of Fez, Box 2202, University S.M. Ben Abdellah Fez, Morocco.
\\ mahdou@hotmail.com }

\author{Mohammed Tamekkante}
\address{Mohammed Tamekkante\\Department of Mathematics, Faculty of Science and Technology of Fez, Box 2202, University S.M. Ben Abdellah Fez, Morocco.
\\ tamekkante@yahoo.fr }



\begin{abstract}
In this paper, we compare the Gorenstein homological dimension of
a ring $R$ and of its trivial ring extension by an module $E$.
\end{abstract}

\maketitle
\section{Introduction}

Throughout this paper, all rings are commutative with identity
element, and all modules are unital.

Let $R$ be a ring, and let $M$ be an $R$-module. As usual we use
$\pd_R(M)$, $\id_R(M)$ and $\fd_R(M)$ to denote, respectively, the
classical projective dimension, injective dimension and flat
dimension of $M$. By $\gldim(R)$ and $\wdim(R)$ we denote,
respectively, the classical global dimension and weak dimension of
R.

Recall that the Gorenstein homological theory starts in the
sixties with Auslander and Bridger \cite{Auslander,Aus bri} over
commutative Noetherian rings and developed, several decades later,
by Enochs, Jenda, Christensen, Holm, Yassemi  and others (see
\cite{Bennis and
Mahdou1,Bennis and Mahdou2,Bennis and Mahdou3,Christensen, Enochs,Enochs2,Holm}).\\

Recently in \cite{Bennis and Mahdou2}, the authors started the
study of global Gorenstein dimensions of rings, which are called,
for a commutative ring $R$, Gorenstein global  projective,
injective, and weak dimensions of $R$, denoted by $GPD(R)$,
$GID(R)$, and $G.wdim(R)$, respectively, and, respectively,
defined as follows:\bigskip

$\begin{array}{cccc}
  1) & GPD(R) & = & sup\{ Gpd_R(M)\mid M$ $R-module\} \\
  2) & GID(R) & = & sup\{ Gid_R(M)\mid M$ $R-module\} \\
  3) & G.wdim(R) & = & sup\{ Gfd_R(M)\mid M$ $R-module\}
\end{array}$
\\

They proved that, for any ring R, $ G.wdim(R)\leq GID(R) = GPD(R)$
(\cite[Theorems 1.1 and Corollary 1.2(1)]{Bennis and Mahdou2}).
So, according to the terminology of the classical theory of
homological dimensions of rings, the common value of $GPD(R)$ and
$GID(R)$ is called
Gorenstein global dimension of $R$, and denoted by $G.gldim(R)$.\\
They also proved that the Gorenstein global and weak dimensions
are refinement of the classical global  and weak dimensions of
rings. That is : $G.gldim(R) \leq gldim(R)$ and $G.wdim(R)\leq
wdim(R)$ with equality if $wdim(R)$ is finite (\cite[Corollary
1.2(2 and 3)]{Bennis and Mahdou2}).\\

Let $R$ be a ring and $E$ an $R$-module. The trivial ring
extension of $R$ by $E$ is the ring $R := A\ltimes E$ whose
underlying group is $A\times E$ with multiplication given by $(r,
e)(r', e') = (rr', re'+r'e)$ (\cite{Fossum, Huckaba, KabbajMahdou,
Kabbaj and Mahdou like conditions}). Over $R\ltimes E$,
 the module $0\times E$ is an ideal. Moreover, the diagonal embedding  $\varphi :R \rightarrow R\bowtie
I$, defined by $\varphi(r)=(r, 0)$, is an injective ring
homomorphism. Hence we have the following short exact sequence of
R-modules:
$$(\ast)\qquad 0\rightarrow R\stackrel{\varphi}\rightarrow R\ltimes
E \stackrel{\psi}\rightarrow E\rightarrow 0$$ where $\psi((r, e))
=e$, for every $(r, e)\in  R \ltimes E$. Notice that this sequence
splits. We also have the short exact sequence of $R\ltimes
E$-modules:
$$(\ast\ast)\qquad 0\longrightarrow 0\times E\stackrel{\iota}\longrightarrow
R\ltimes  E\stackrel{\varepsilon}\longrightarrow R \longrightarrow
0$$ where $i$ is the injection and $\varepsilon(r,e)=r$. Note that
$R$ is an $R\ltimes  E$-module via the map ring $\varepsilon$
(explicitly for all $r,r'\in R$ and $e\in E$,
$(r,e).r'=\varepsilon(r,e)r'=rr'$).\\
Contrarily to $(\ast)$ this sequence never splits as shown by the
following result:
\begin{proposition}
Let $R$ be a ring and $E\neq 0$ an $R$-module. Then, $R$ is never
projective as an $R\ltimes E$-module.
\end{proposition}
\begin{proof}
Consider the short exact sequence of $R\ltimes  E$-modules:
$$(\ast\ast)\qquad 0\longrightarrow 0\times E\stackrel{\iota}\longrightarrow
R\ltimes  E\stackrel{\varepsilon}\longrightarrow R \longrightarrow
0$$ where $i$ is the injection and $\varepsilon(r,e)=r$. It is
clear that $R$ is projective if, and only if, $(\ast\ast)$ splits.
Hence, there is an $R\ltimes E$-morphism $\pi: R\rightarrow
R\ltimes E$ such that $\varepsilon\circ \pi=id(R)$. Set
$\pi(1)=(r,e_0)$. Thus, $1=\varepsilon\circ
\pi(1)=\varepsilon(r,e_0)=r$. Hence, for an arbitrary $r\in R$ and
any $e\in E$, we have
$\pi(r)=\pi((r,e).1))=(r,e)(1,e_0)=(r,re_0+e)$. But that is
impossible since $\pi$ must be well defined and $E\neq 0$.
\end{proof}
More general, If $E$ is a flat $R$-module, from \cite[Corollary
4.7]{Fossum}, we conclude:
\begin{lemma}\label{lemma0}
Let $R$ be a ring and $E$ a flat $R$-module. Then, $fd_{R\ltimes
E}(R)\leq n$ if, and only if, $E^n:=\underbrace{E\otimes E\otimes
...\otimes E}_n=0$.
\end{lemma}

To give examples of Lemma \ref{lemma0}, we have to think about
rings which contain an idempotent element (i.e; $a\in R$ such that
$a^n=0$ for a positive integer $n$).\\

The homological behavior and structure of the $R\ltimes E$-module
$R$ has an importance counterpart in the determination of the
Gorenstein and classical dimensions of the ring $R\ltimes E$.
Recall that (see \cite[Proposition 2.6]{Bennis and Mahdou2})
$$R\;  is\;  quasi-Frobenius \quad\Longleftrightarrow \quad G.gldim(R)=0 $$
 Adding the Noetherian condition to \cite[Corollary 4.36]{Fossum}
 we obtain the next corollary:
 \begin{corollary}
 Let $R$ be a Noetherian ring and $E$ a finitely generated
 $R$-module. Then, $R\ltimes E$ is quasi-Frobenius if, and only
 if, the following conditions hold:
 \begin{enumerate}
    \item $E$ and $Ann_R(E)$ are injective $R$-modules,
    \item  The naturel map $R\rightarrow Hom_R(E,E)$ is an
    epimorphism, and
    \item $Hom_R(E,Ann_R(E))=0$.
 \end{enumerate}
 \end{corollary}

 In \cite{MahoduOuarghi}, the authors study the Gorenstein
 dimension in trivial ring extensions. Namely they proved that for an
 $R$-module $E$ with finite flat dimension such that $G.gldim(R)<\infty$, we have
 $G.gldim(R)\leq G.gldim(R\ltimes E)+fd_R(E)$
 (\cite[Theorem 2.4]{MahoduOuarghi}). Moreover, in
 \cite{MahoduOuarghi} we find examples of trivial ring extensions of
 rings with infinite Gorenstein global dimension (see
 \cite[Theorems 3.2 and 3.4]{MahoduOuarghi}).\\

In this paper we need the condition $Tor^n_{R\ltimes E}(R,M)=0$
for any positive integer $n$ and all $R\ltimes E$-module with
finite projective dimensions. To give an example of this situation
we take   $E=xR$  where $x$ is a nonzero divisor. We have the
short exact sequences of $R\ltimes xR$-modules:
$$(1)\quad 0\longrightarrow 0\times xR\stackrel{\iota}\longrightarrow R\ltimes
xR\stackrel{\psi}\longrightarrow R\longrightarrow 0$$
$$(2)\quad 0\longrightarrow R \stackrel{\mu}\longrightarrow R\ltimes xR \stackrel{\nu}\longrightarrow
0\times xR \longrightarrow 0$$  where $\iota$ is the injection,
$\psi((r,xr')=r$, $\mu(r)=(0,xr)$ and $\nu(r,xr')=(0,xr)$.  Then,
from $(1)$ and $(2)$, for every $R\ltimes E$-module with finite
projective dimension we have $Tor^i_{R\ltimes xR}(R,M)=0$.

\begin{discussion}[\textbf{Modules over $R\ltimes E$}]\label{discussion1}
Via  the ring map $\varepsilon: R\ltimes E\rightarrow R$ defined
by $\varepsilon(r,e)=r$, we can give every $R\ltimes E$-module $M$
a structure of $R$-module  (by setting $r.m:=(r,0)m$ for every
$r\in R$). Moreover, we can consider  the $R$-morphism (which
depend only to the modulation  of $M$ over $R\ltimes  E$); $\rho:
E\otimes M \rightarrow M$ defined by $\rho(e\otimes m)=(0,e)m$
(see that $\rho$ is well defined by the universal propriety of
tensor product).  This $R$-morphism satisfying the condition
$$(\mathcal{H}):\qquad
  \rho(e\otimes
\rho(e'\otimes m))=0\qquad
 for\; e,e'\in E\; and\; m\in M$$ Conversely, given a pair
$(M,\rho)$ where $M$ is an $R$-module and $\rho:
E\otimes_RM\rightarrow M$ is an $R$-morphism which satisfied the
condition $(\mathcal{H})$. We can give $M$ an $R\ltimes E$-module
structure via $\rho$. Namely for $e\in E$, $r\in R$ and $m\in M$
$$(r,e).m:=rm+\rho(e\otimes m)$$ (to see that the condition $(\mathcal{H})$ guaranties the modulation we have just to try to prove that
$(r,e)[(r',e').m]=[(r,e)(r',e')].m$). These two constructions are
inverse of each other. Hence,   we can identified
an $R\bowtie E$-module M to a pair $(M,\rho)$ where $\rho$ satisfies the condition $(\mathcal{H})$. \\
A revealing example is to examine how $R\ltimes E$ is identifying
to a pair $(R\ltimes, \rho)$. So, as above we define $\rho$ as
$\rho(e\otimes(r,e'))=(0,e)(r,e')=(0,re)$.\\

\end{discussion}
Recall some  well-know results.
\begin{proposition}\label{facts}
Let $R$ be a ring and $E$ be an $R$-module. Then, for any
$R$-module $M$ we have:
\begin{enumerate}
  \item $pd_R(M)\leq pd_{R\ltimes E}(M)$,
  \item $id_R(M)\leq id_{R\ltimes E}(M)$, and
  \item $fd_R(M)\leq fd_{R\ltimes E}(M)$.
\end{enumerate}
\end{proposition}
\begin{proof}
$(1)$ and $(2)$ are the particular cases of \cite[Lemmas 4.1 and
4.2]{Fossum}.\\
$(3)$ If $fd_{R\ltimes E}(M)<\infty$, then by \cite[Lemma 3.51 and
Theorem 3.52]{Rotman}, we have $id_{R\ltimes
E}(Hom_{\mathbb{Z}}(M,\mathbb{Q}/\mathbb{Z}))=fd_{R\ltimes E}(M)$.
But the $R\ltimes E$-modulation over
$Hom_{\mathbb{Z}}(M,\mathbb{Q}/\mathbb{Z})$ is defined by, for
every $(r,e)\in R\ltimes E$ and $f\in
Hom_{\mathbb{Z}}(M,\mathbb{Q}/\mathbb{Z})$, $(r,e).f:M\rightarrow
\mathbb{Q}/\mathbb{Z}$ such that for any $m\in M$,
$$((r,e).f)(m)=f((r,e).m)=f(rm)=rf(m)$$ Thus, $(r,e).f=rf$. Hence, by
$(2)$ obove, $id_R(Hom_{\mathbb{Z}}(M,\mathbb{Q}/\mathbb{Z}))\leq
id_{R\ltimes E}(Hom_{\mathbb{Z}}(M,\mathbb{Q}/\mathbb{Z}))$.
Consequently, we have:
$$fd_R(M)=id_R(Hom_{\mathbb{Z}}(M,\mathbb{Q}/\mathbb{Z}))\leq id_{R\ltimes E}(Hom_{\mathbb{Z}}(M,\mathbb{Q}/\mathbb{Z}))=fd_{R\ltimes E}(M).$$
\end{proof}

\section{Main results}

The aim  of this paper is to  give a Gorenstein version of
Proposition \ref{facts}.

\begin{theorem}\label{first main result}
Let $R$ be a ring and $E$ an $R$ -module such that $pd_{R\ltimes
E}(R)<\infty$. Then, for any $R$-module $M$ we have $Gpd_R(M)\leq
Gpd_{R\ltimes E}(M)$. Consequently, $G.gldim(R)\leq
G.gldim(R\ltimes E).$
\end{theorem}
To prove this Theorem we involve several Lemmas.
\begin{lemma}\label{lemma Gorenstein}
Let $R$ be a ring and $E$ an $R$-module such that $pd_{R\ltimes
E}(R)<\infty$. If $M$  is a Gorenstein projective
 $R\ltimes E$-module then $M\otimes_{R\ltimes E}R$ is a Gorenstein
 projective $R$-module. Moreover, if $Tor^i_{R\ltimes
 E}(M,R)=0$ for all $i>0$ then $Gpd_R(M\otimes_{R\ltimes E}R)\leq
 Gpd_{R\ltimes E}(M)$.
 \end{lemma}
\begin{proof}
Note in first that for every Gorenstein projective $R\ltimes
E$-module, we have $Tor_{R\ltimes E}(M,R)=0$. Indeed, we can pick
an exact sequence of $R\ltimes E$-modules:
$$0\rightarrow M \rightarrow P_1\rightarrow P_2\rightarrow
...\rightarrow P_n\rightarrow G\rightarrow 0$$ where all $P_i$ are
projective and for any integer $n>0$ (in particular for
$n:=pd_{R\ltimes E}(R)$). Hence, $Tor_{R\ltimes
E}(M,R)=Tor_{R\ltimes E}^{n+1}(G,R)=0$. Recall also that a
Gorenstein projective module is an image of
a morphism in a complete projective resolution.\\
Let $M$ be an arbitrary $R\ltimes E$-module and consider a
complete projective resolution of $R\ltimes E$-modules:
$$\mathbf{P}:...\rightarrow P_1\rightarrow
P_0\rightarrow P^0\rightarrow P^1\rightarrow...$$ such that
$M=Im(P_0\rightarrow P^0)$ (\cite[Definition 2.1]{Holm}). By the
reason above, the operator $-\otimes_{R\ltimes E}R$ leaves
$\mathbf{P}$ exact. Then, we obtain an exact sequence of
$R$-modules:
$$\mathbf{P}\otimes_{R\ltimes E}R:...\rightarrow P_1\otimes_{R\ltimes E}R\rightarrow
P_0\otimes_{R\ltimes E}R\rightarrow P^0\otimes_{R\ltimes
E}R\rightarrow P^1\otimes_{R\ltimes E}R\rightarrow...$$ On the
other hand, for each projective $R$-module $Q$ we have
$pd_{R\ltimes E}(Q)\leq pd_{R\ltimes E}(R)<\infty$. Thus,
$Hom_R(\mathbf{P}\otimes_{R\ltimes E}R,Q)\cong Hom_{R\ltimes
E}(\mathbf{P},Q)$ is exact (\cite[Proposition 2.3]{Holm}).
Consequently, $\mathbf{P}\otimes_{R\ltimes E}R$ is a complete
projective resolution of $R$-modules. So, $M\otimes_{R\ltimes
E}R\cong Im((P_0\rightarrow P^0)\otimes 1_R)$
is a Gorenstein projective $R$-module, as desired.\\
Now let $M$ be an $R\ltimes E$-module with finite Gorenstein
projective dimension equal to $n$ such that $Tor^i_{R\ltimes
E}(M,R)=0$ for all $i>0$. The desired result follows by applying
the functor $-\otimes_{R\ltimes E}R$ to an $n$-step Gorenstein
projective resolution of $M$ over $R\ltimes E$.
\end{proof}
\begin{lemma}\label{lemma suite}
Let $0\rightarrow N\rightarrow N' \rightarrow N'' \rightarrow 0$
be an exact sequence of $R$-modules. Then,
   $Gpd_R(N'')\leq max\{Gpd_R(N'),Gpd_R(N)+1\}$ with equality
  if $Gpd_R(N')\neq Gpd_R(N)$.
\end{lemma}
\begin{proof}
Using \cite[Theorems 2.20 and 2.24]{Holm} the argument is
analogous to the one of \cite[Corollary 2, p. 135]{bourbaki}.
\end{proof}
\begin{lemma}\label{proposition G-projective} Let $R$ be a ring and
$E$ an $R$-module such that $pd_{R\ltimes E}(R)<\infty$. Let $B$
and $D$ a couple of $R$ modules and $\rho: E\otimes_R (B\oplus
D)\rightarrow B\oplus D$ which satisfies the condition
$(\mathcal{H})$ (see Discussion \ref{discussion1}) and such that
$Im(\rho)\subseteq 0\oplus D$. With the identification of
Discussion \ref{discussion1}, we have  $Gpd_R(B)\leq Gpd_{R\ltimes
E}((B\oplus D,\rho))$.
\end{lemma}
\begin{proof}
Recall that the $R\ltimes E$-modulation over $(B\oplus
 D,\rho)$ is given by setting:
 $$(r,e).(b,d):=r(b,d)+\rho(e\otimes (b,d))$$
 (see Discussion \ref{discussion1}).\\
In first, we assume that $(B\oplus D,\rho)$ is a Gorenstein
projective $R\ltimes E$-module and we claim that $B$ is a
Gorenstein projective $R$-module. Seeing that $Im(\rho)=(0\times
I)(B\oplus D)$ and since $R\cong R\ltimes E/(0\times E)$, it is
clear that $(B\oplus D)/Im(\rho)\cong (B\oplus
D,\rho)\otimes_{R\ltimes E}R$ is a Gorenstein projective
$R$-module (by Lemma \ref{lemma Gorenstein}). Now, consider the
$R$-morphisms: $(B\oplus D)/Im(\rho)\stackrel{\delta}\rightarrow
B$ and $B\stackrel{\delta'}\rightarrow (B\oplus D)/Im(\rho)$
defined by $\delta(\overline{(b,d)})=b$ and
$\delta'(b)=\overline{(b,0)}$. We can see easily that $\delta$ is
well defined. Indeed, if $\overline{(b,d)}=\overline{(b',d')}$
then $(b-b',d-d')\in Im(\rho )$ and so, $b-b'=0$ (since
$Im(\rho)\subseteq 0\oplus D$). Also, we can check  that
$\delta\circ\delta'=id(B)$. Then, $B$ is a direct summand of
$(B\oplus D)/Im(\rho)$. Hence,  $B$ is  a Gorenstein projective
$R$-module (by \cite[Theorem 2.5]{Holm}). Therefore, we assume
$0<n:=pd_{R\ltimes E}((B\oplus D,\rho))$ and we
proceed by induction on $n$.\\
Inspecting the proof of \cite[Lemma 4.1]{Fossum} we can construct
a short exact sequence of $R\ltimes E$-modules with the form
$$0\longrightarrow (K\oplus L,\phi)\longrightarrow Q \longrightarrow
(B\oplus D,\rho)\longrightarrow 0$$ where $Q$ is projective and
$Im(\phi)\subseteq 0\oplus L(:=L)$. Hence, by the hypothesis
induction and Lemma \ref{lemma suite}, we conclude that:
$$Gpd_{R\bowtie I}(B\oplus D,\rho)=1+Gpd_{R\bowtie I}(K\oplus
L,\phi)\geq 1+Gpd_R(K)\geq Gpd_R(B)$$
\end{proof}
\begin{proof}[Proof of Theorem \ref{first main result}]
Recall that the modulation of $R\bowtie I$ over the $R$-module $M$
is defined via the ring map $R\bowtie I\rightarrow R$ defined by
$(r,r+i)\mapsto r$. Explicitly, we have for all $m\in M$,
$(r,r+i).m=rm$. So, we can identify this $R\bowtie I$-module with
the $R\bowtie I$-module $(M,\rho)$ with $\rho:I\otimes
M\rightarrow M$ is the zero $R$-morphism. Thus, by Lemma
\ref{proposition G-projective}, $pd_R(M)\leq pd_{R\bowtie I}(M)$,
as desired.
\end{proof}

\begin{remark}
Notice that the hypothesis of Theorem \ref{first main result} is
sufficient but not necessary. A simple example to see that is by
considering the ring $R\ltimes R$ where $R$ is coherent. Using
\cite[Theorem
 1.4.5]{Glaz} we can prove that  $fd_{R\ltimes R}(R)=\infty$. But,
 $G.gldim(R\ltimes R)=G.gldim(R)$ (\cite[Proposition 2.5]{Bennis and
 Mahdou4}). In \cite[Proposition 3.11 and Corollary 5.5]{Fossum}, the authors give an other example of our remark.
 Namely, if $E$ is a finitely generated projective module over a Noetherian ring $R$ then,
$$ R\ltimes E \; is\; n-Gorenstein \quad \Longrightarrow \quad R\;
is \; n-Gorenstein$$

 Recall that $R$ is called  $n$-Gorenstein if it is
 Noetherian with $id_{R}(R)\leq n$ and note that if  $R$ is a Noetherian ring then
 $G.gldim(R)\leq n \Leftrightarrow R$ is $n$-Gorenstein ( for
 $\Rightarrow$ see \cite{Enochs and Xu} and for $\Leftarrow$ use
 \cite[Theorem 2.20]{Holm}).
\end{remark}

\begin{proposition}\label{main result 2}
Let $R$ be a ring and $E$ an $R$ module such that
$G.gldim(R\ltimes E)<\infty$. Suppose that $Tor^i_{R\ltimes
E}(M,R)=0$ for all $i>n$ and every $R\ltimes E$-module $M$ with
finite projective dimension. Then, $G.gldim(R\ltimes E)\leq
G.gldim(R)+n$.
\end{proposition}

To prove this Proposition, we need the following Lemma.

\begin{lemma}\label{lemma2}
Let $R$ be a ring with  finite Gorenstein projective dimension,
then, for a positive integer $n$, the following statements are
equivalent:
    \begin{enumerate}
        \item $G.gldim(R)\leq n$;
        \item $pd(I)\leq n$ for every injective module $I$.
    \end{enumerate}
    \end{lemma}
    \begin{proof} Note that
    $G.gldim(R)=sup \{ Gid(M)|M\; an\; R-module\}$  (by \cite[Theorem
    1.1]{Bennis and Mahdou2}). Thus, using \cite[Theorem 2.22]{Holm},
    $G.gldim(R)\leq n \Leftrightarrow Ext^i(I,M)=0$ for each $i>n$ and
    for any injective module $I$ and each module $M$. Thus,
    $G.gldim(R)\leq n \Leftrightarrow pd(I)\leq n$ for each injective
    module $M$, as desired.
    \end{proof}
\begin{proof}[Proof of Proposition \ref{main result 2}]
We may assume that $m:=G.gldim(R)$ and $n$ are finite. Otherwise,
the result is obvious. Let $I$ be an arbitrary injective $R\ltimes
E$-module. Since $G.gldim(R\ltimes E)<\infty$, we have
$pd_{R\ltimes E}(I)<\infty$ (by Lemma \ref{lemma2}). For such
module pick an $n$-step projective resolution as follows:
$$0\rightarrow K \rightarrow P_n\rightarrow ... \rightarrow
P_1\rightarrow I \rightarrow 0$$ Hence, $Tor^i_{R\ltimes
E}(K,R)=0$ for all $i>0$. Thus, using \cite[Theorem 4.9]{Fossum},
$pd_R(K\otimes_{R\ltimes E}R)=pd_{R\ltimes E}(K)<\infty$. Then,
$pd_{R\ltimes E}(K)=pd_R(K\otimes_{R\ltimes E}R)\leq m$ (by
\cite[Corollary 2.7]{Bennis and Mahdou2}). Consequently,
$pd_{R\ltimes E}(I)\leq G.gldim(R)+n$. Thus, from Lemma
\ref{lemma2}, we obtain the desired result.
\end{proof}

\begin{corollary}
Let $R$ be a ring and $E$ a non-zero cyclic $R$ module such that
$G.gldim(R\ltimes E)<\infty$. Then, $G.gldim(R\ltimes E)\leq
G.gldim(R)$.
\end{corollary}
\begin{proof}
Inspecting the proof of \cite[Theorem 2.28]{Fossum} we see that
for a cyclic $R$-module $E$ we have: $Tor^i_{R\ltimes E}(M,R)=0$
for all $i>0$ and each $R\ltimes E$-module with finite projective
dimension $M$. Thus, the desired result follows directly from
Proposition \ref{main result 2}.
\end{proof}

 Now we give our second main result
in this paper.

\begin{theorem}\label{main result3}
Let $R$ be a ring and $E$ an $R$-module such that $R\ltimes E$ is
coherent and such that $fd_{R\ltimes E}(R)<\infty$. Then, for any
$R$-module $M$ we have $Gfd_R(M)\leq Gfd_{R\ltimes E}(M)$.
Consequently, $G.wdim(R)\leq G.wdim(R\ltimes E)$.
\end{theorem}

First we have to recall that in \cite[Theorem 4.4.4]{Glaz}, Glaz
gives the necessary and sufficient condition under $R$ and $E$ to
obtain the coherence of $R\ltimes E$ and  make sure that if
$R\ltimes E$ is coherent, so is $R$.

\begin{lemma}\label{Gorenstein flat}
Let $R$ be a ring and $E$ an $R$-module such that $fd_{R\ltimes
E}(R)<\infty$. If $M$ is a Gorenstein flat
 $R\ltimes E$-module then $M\otimes_{R\ltimes E}R$ is a Gorenstein
 flat  $R$-module.
 \end{lemma}

 \begin{proof} Let $M$ be a Gorenstein flat $R\ltimes
 E$-module and consider  a  complete flat resolution of $R\ltimes E$-modules:

$$ \mathbf{F}: ...\rightarrow F_1\rightarrow F_0 \rightarrow F^0\rightarrow
F^1\rightarrow...$$ such that $M=Im(F_0\rightarrow F^0)$
(\cite[Definition 3.1]{Holm}).
 By the same reason that in the proof of Lemma \ref{lemma
 Gorenstein},
 the operator $-\otimes_{R\ltimes E}R$ leaves $\mathbf{F}$ exact since  $fd_{R\ltimes E}(R)<\infty$. So, we obtain  the exact flat resolution of
 $R$-modules:
 $$\mathbf{F}\otimes_{R\ltimes E}R: ...\rightarrow F_1\otimes_{R\ltimes E}R\rightarrow F_0\otimes_{R\ltimes E}R \rightarrow F^0\otimes_{R\ltimes E}R\rightarrow F^1\otimes_{R\ltimes E}R\rightarrow...$$
 Now  let $I$ be an injective $R$-module, $N$
an arbitrary $R\ltimes E$-module and set $fd_{R\ltimes E}(R)=n$.
Pick an $n$-step projective resolution of $N$ over $R\ltimes E$ as
follows:
$$0\rightarrow N'\rightarrow P_n\rightarrow...\rightarrow
P_1\rightarrow N\rightarrow 0$$ Clearly, $Tor_{R\ltimes
E}(N',R)=Tor^{n+1}_{R\ltimes E}(N,R)=0$. Thus, from
\cite[Proposition 4.1.3]{Cartan}, we have $Ext_{R\ltimes
E}(N',I)\cong Ext_R(N'\otimes_{R\ltimes E}R,I)=0$. Therefore,
$Ext^{n+1}_{R\ltimes E}(N,I)=Ext_{R\ltimes E}(N',I)=0$.
Consequently, $id_{R\ltimes E}(I)\leq n<\infty$. Then, the complex
$\mathbf{F}\otimes_{R\ltimes
E}R\otimes_RI\cong\mathbf{F}\otimes_{R\ltimes E}I$
 is exact (direct consequence of \cite[Theorem 3.14]{Holm}) and so $\mathbf{F}\otimes_{R\ltimes
E}R$ is a complete flat resolution of $R$-modules. Therefore,
$M\otimes_{R\ltimes E}R=Im(F_0\otimes_{R\ltimes E}R \rightarrow
F^0\otimes_{R\ltimes E}R)$ is a Gorenstein flat module.
\end{proof}

Using \cite[Proposition 3.11]{Holm} and the injective version of
Lemma \ref{lemma suite} we get the following Lemma:

\begin{lemma}\label{lemma suite2}
Let $0\rightarrow N\rightarrow N' \rightarrow N'' \rightarrow 0$
be an exact sequence of modules over a coherent ring $R$. Then:
$Gfd_R(N'')\leq max\{Gfd_R(N'),Gfd_R(N)+1\}$ with equality
  if $Gfd_R(N')\neq Gfd_R(N)$.
\end{lemma}
\begin{proof}[Proof of Theorem \ref{main result3}] Recall that $R$
is also coherent (by \cite[Theorem 4.4.4]{Glaz}). Similarly that
in the proof of Lemma \ref{proposition G-projective}; by replacing
Lemma \ref{lemma Gorenstein}, \cite[Theorem 2.5]{Holm} and Lemma
\ref{lemma suite} by Lemma \ref{Gorenstein flat},
\cite[Proposition 3.13]{Holm} and Lemma \ref{lemma suite2}
respectively, we prove that: if  $B$ and $D$ are a couple of $R$
modules and $\rho: E\otimes_R (B\oplus D)\rightarrow B\oplus D$
which satisfies the condition $(\mathcal{H})$ (see Discussion
\ref{discussion1}) and such that $Im(\rho)\subseteq 0\oplus D$,
then $Gfd_R(B)\leq Gfd_{R\ltimes E}((B\oplus D,\rho))$.
Consequently, as in the proof of Theorem \ref{first main result},
we deduce that for any $R$-module $M$ we have: $Gfd_R(M)\leq
Gfd_{R\ltimes E}(M)$, as desired.
\end{proof}

\bibliographystyle{amsplain}

\end{document}